\documentclass{birkjour}
\usepackage[T1]{fontenc}
\usepackage[latin9]{inputenc}
\usepackage{amsthm}
\usepackage{amstext}
\usepackage{amssymb}

\makeatletter
\numberwithin{equation}{section}
\numberwithin{figure}{section}
\theoremstyle{plain}
\newtheorem{thm}{\protect\theoremname}
  \theoremstyle{plain}
  \newtheorem{lem}[thm]{\protect\lemmaname}
  \theoremstyle{plain}
  \newtheorem{cor}[thm]{\protect\corollaryname}
  \theoremstyle{plain}
  \newtheorem{prop}[thm]{\protect\propositionname}

\usepackage{amscd}

\makeatother

\usepackage[english]{babel}
  \providecommand{\corollaryname}{Corollary}
  \providecommand{\lemmaname}{Lemma}
  \providecommand{\propositionname}{Proposition}
\providecommand{\theoremname}{Theorem}

\begin{document}
\global\long\def\cc{\mathbb{C}}

\global\long\def\pm{\mathbb{P}^{m}}

\title{Reduction of a family of ideals}

\author{Tomasz Rodak}
\address{University of \L{}ód\'{z}\\
Faculty of Mathematics and Computer Science\\
S. Banacha 22, 90-238 \L{}ód\'{z}, Poland}
\email{rodakt@math.uni.lodz.pl}
\thanks{This research was partially supported by the Polish OPUS Grant No
2012/07/B/ST1/03293}

\begin{abstract}
In the paper we prove that there exists a simultaneous reduction of
one-parameter family of $\mathfrak{m}_{n}$-primary ideals in the
ring of germs of holomorphic functions. As a corollary we generalize
the result of A. P\l{}oski \cite{ploski} on the semicontinuity of
the \L{}ojasiewicz exponent in a multiplicity-constant deformation.
\end{abstract}

\maketitle

\section{Introduction}

Let $R$ be a ring and $I$ an ideal. We say that an ideal $J$ is
a \emph{reduction of }$I$ if it satisfies the following condition:
\[
J\subset I,\quad\text{and for some}\quad r>0\quad\text{we have}\quad I^{r+1}=JI^{r}.
\]

The notion of reduction is closely related to the notions of \emph{Hilbert-Samuel
multiplicity} and \emph{integral closure} of an ideal.

Recall that if $\left(R,\mathfrak{m}\right)$ is a Noetherian local
ring of dimension $n$ and $I$ is an $\mathfrak{m}$-primary ideal of
$R$, then the \emph{Hilbert-Samuel multiplicity }of $I$ is given
by the formula
\[
e(I)=n!\lim_{k\to\infty}\frac{\mathrm{length}_{R}R/I^{k}}{k^{n}}.
\]

For the multiplicity theory in local rings see for example \cite{matsumura}
or \cite{huneke}.

Let $I$ be an ideal in a ring $R$. An element $x\in R$ is said
to be \emph{integral over }$I$ if there exists an integer $n$ and
elements $a_{k}\in I^{k}$, $k=1,\ldots,n$, such that 
\[
x^{n}+a_{1}x^{n-1}+\cdots+a_{n}=0.
\]

The set of all elements of $R$ that are integral over $I$ is called
the \emph{integral closure of }$I$, and is denoted $\overline{I}$.
If $I=\overline{I}$ then $I$ is called \emph{integrally closed}.
It is well known that $\overline{I}$ is an ideal. 

The relationship between the above notions is given in the following
Theorem due to D. Rees:
\begin{thm}[{Rees, \cite[Cor. 1.2.5, Thm. 11.3.1]{huneke}}]
\label{thm:Rees}Let $(R,\mathfrak{m})$ be a formally equidimensional
Noetherian local ring and let $J\subset I$ be two $\mathfrak{m}$-primary
ideals. Then the following conditions are equivalent:
\begin{enumerate}
\item $J$ is a reduction of $I$;
\item $e\left(I\right)=e\left(J\right)$;
\item $\overline{I}=\overline{J}$.
\end{enumerate}
\end{thm}
It is an important fact that a reduction of an ideal is often generated
by a system of parameters. More precisely we have
\begin{thm}[{\cite[Theorem 14.14]{matsumura}}]
\label{thm:rzutowanie}Let $(R,\mathfrak{m})$ be a $d$-dimensional
Noetherian local ring, and suppose that $k=R/\mathfrak{m}$ is an infinite
field; let $I=(u_{1},\ldots,u_{s})$ be an $\mathfrak{m}$-primary ideal.
Then there exist a finite number of polynomials $D_{\alpha}\in k[Z_{ij};1\leqslant i\leqslant d,1\leqslant j\leqslant s]$,
$1\leqslant\alpha\leqslant\nu$ such that if $y_{i}=\sum a_{ij}u_{j}$,
$i=1,\ldots,d$ and at least one of $D_{\alpha}(\overline{a}_{ij};1\leqslant i\leqslant d,1\leqslant j\leqslant s)\ne0$,
then the ideal $(y_{1},\ldots,y_{d})R$ is a reduction of $I$ and
$\{y_{1},\ldots,y_{d}\}$ is a system of parameters of $R$.
\end{thm}
Let $\left(\mathcal{O}_{n},\mathfrak{m}_{n}\right)$ be the ring of
germs of holomorphic functions $\left(\mathbb{C}^{n},0\right)\to\mathbb{C}$.
The aim of this note is to prove the following:
\begin{thm}
\label{thm:Main}Let $F=F_{t}(x)=F(x,t)\colon(\cc^{n}\times\cc,0)\to(\cc^{m},0)$
be a holomorphic map. Assume that $(F_{t})\mathcal{O}_{n}$ is an
$\mathfrak{m}_{n}$-primary ideal for all $t$. Then there exists
a complex linear map $\pi\colon\mathbb{C}^{m}\to\mathbb{C}^{n}$ such
that for all $t$ the ideal $\left(\pi\circ F_{t}\right)\mathcal{O}_{n}$
is a reduction of $\left(F_{t}\right)\mathcal{O}_{n}$.
\end{thm}
In the next section we get as a corollary that if the above family
$\left(F_{t}\right)\mathcal{O}_{n}$ is of constant multiplicity then
the \L{}ojasiewicz exponent in this family is a lower semicontinuos
function of $t$. A. P\l{}oski proved this result under additional
restriction $m=n$ but with space of parameters of arbitrary dimension.

The proof of Theorem \ref{thm:Main} is based on some geometric property
of Hilbert-Samuel multiplicity, given in section 3.

\section{Semicontinuity of the \L{}ojasiewicz exponent}

Let $(R,\mathfrak{m})$ be a local ring and let $I$ be an $\mathfrak{m}$-primary
ideal. By the \emph{\L{}ojasiewicz exponent }$\mathcal{L}(I)$ of
$I$ we define the infimum of 
\[
\left\{ \frac{p}{q}:\mathfrak{m}^{p}\subset\overline{I^{q}}\right\} .
\]
It was proved in \cite{ljt1974} that if $F\colon\left(\cc^{n},0\right)\to\left(\cc^{m},0\right)$
is a holomorphic map with an isolated zero at the origin and $I:=\left(F\right)\mathcal{O}_{n}$,
then $\mathcal{L}\left(I\right)$ is an optimal exponent $\nu$ in
the inequality
\[
\left|F\left(x\right)\right|\geqslant C\left|x\right|^{\nu},
\]
where $C$ is some positive constant and $x$ runs through sufficiently
small neighbourhood of $0\in\cc^{n}$. 
\begin{lem}
\label{lem:wykladnik - redukcja}Let $\left(R,\mathfrak{m}\right)$
be a Noetherian local ring. If $I$ is an $\mathfrak{m}$-primary
ideal of $R$ and $J$ is a reduction of $I$ then $\mathcal{L}\left(I\right)=\mathcal{L}\left(J\right)$. \end{lem}
\begin{proof}
Obviously $\mathcal{L}(I)\leqslant\mathcal{L}(J)$. Assume that $\mathfrak{m}^{p}\subset\overline{I^{q}}$.
Since $J$ is a reduction of $I$, then also $J^{q}$ is a reduction
of $I^{q}$ \cite[Prop. 8.1.5]{huneke}. Thus $\overline{J^{q}}=\overline{I^{q}}$
by Theorem \ref{thm:Rees}, which gives $\mathfrak{m}^{p}\subset\overline{J^{q}}$.
This proves the inequality $\mathcal{L}(J)\leqslant\mathcal{L}(I)$
and ends the proof.\end{proof}
\begin{cor}[A. P\l{}oski for $m=n$, \cite{ploski}]
Let $F\colon\left(\cc^{n}\times\cc,0\right)\to\left(\cc^{m},0\right)$
be a holomorphic map. Put $I_{t}:=\left(F_{t}\right)\mathcal{O}_{n}$.
If the function $t\mapsto e\left(I_{t}\right)$ is constant and finite
then the function $t\mapsto\mathcal{L}\left(I_{t}\right)$ is lower
semicontinuos.\end{cor}
\begin{proof}
By Theorem \ref{thm:Main} there exists a linear map $\pi\colon\cc^{m}\to\cc^{n}$
such that $J_{t}:=\left(\pi\circ F_{t}\right)\mathcal{O}_{n}$ is
a reduction of $I_{t}$ for all $t$. Thus $\mathcal{L}\left(J_{t}\right)=\mathcal{L}\left(I_{t}\right)$
and $e\left(J_{t}\right)=e\left(I_{t}\right)$ by Theorem \ref{thm:Rees}
and Lemma \ref{lem:wykladnik - redukcja}. Consequently $t\mapsto e\left(J_{t}\right)$
is constant and finite and the assertion follows from the case $m=n$
proved by A. P\l{}oski.
\end{proof}

\section{Improper intersection multiplicity}

Let $I$ be an $\mathfrak{m}_{n}$-primary ideal of $\mathcal{O}_{n}$
and let $f_{1},\ldots,f_{m}$ be its generators. We put $f=\left(f_{1},\ldots,f_{m}\right)$.
It is well known that if $m=n$ then 
\[
e(I)=\dim_{\cc}\mathcal{O}_{n}/I.
\]
On the other hand, if $m>n$ then we may define so-called \emph{improper
intersection multiplicity }$i_{0}\left(I\right)$ of $I$ as the improper
intersection multiplicity $i(\mathrm{graph}f\cdot(\cc^{n}\times\{0\});(0,0))$
of $\mathrm{graph}f$ and $\cc^{n}\times\{0\}$ at the point $(0,0)\in\cc^{n}\times\cc^{m}$
(see \cite{atw}). 

Let $C_{f}$ be the (Whitney) tangent cone of the germ of the image
of $f$ at the origin. The following observation is due to S. Spodzieja.
\begin{thm}[\cite{spodzieja-mult}]
\label{thm:sp}The number $i_{0}\left(I\right)$ is well defined.
Moreover, if $\pi\colon\cc^{m}\to\cc^{l}$ is a linear map such that
$\ker\pi\cap C_{f}=\{0\}$, then the ideal $J$ generated by $\pi\circ f$
is $\mathfrak{m}_{n}$-primary and we have $i_{0}(I)=i_{0}(J)$. If additionally
$l=n$ then $i_{0}(I)=e(J)$.\end{thm}
\begin{cor}
\label{prop:i0=00003De}If $I$ is an $\mathfrak{m}_{n}$-primary ideal
in $\mathcal{O}_{n}$, then $i_{0}(I)=e(I)$. \end{cor}
\begin{proof}
Let $I=(f_{1},\ldots,f_{m})\mathcal{O}_{n}$. By Theorems \ref{thm:rzutowanie}
and \ref{thm:sp} there exists linear combinations $g_{i}=\sum a_{ij}f_{j}$,
$i=1,\ldots,n$ such that $J=(g_{1},\ldots,g_{n})\mathcal{O}_{n}$
is a reduction of $I$, $\{g_{1},\ldots,g_{n}\}$ is a system of parameters
of $\mathcal{O}_{n}$ and $i_{0}(I)=i_{0}(J)=e(J)$. From Theorem
\ref{thm:Rees} we get $e(I)=e(J)$. This ends the proof.\end{proof}
\begin{cor}
\label{cor:cone-reduction}If $\pi\colon\cc^{m}\to\cc^{l}$ is a linear
map such that $\ker\pi\cap C_{f}=\{0\}$, then the ideal $J$ generated
by $\pi\circ f$ is a reduction of $I$.\end{cor}
\begin{proof}
We have $J\subset I$ and $e\left(J\right)=e\left(I\right)$. This
and Theorem \ref{thm:Rees} give the assertion.
\end{proof}

\section{Elementary blowing-up}

Here we recall the notion of an elementary blowing-up after \cite{lojasiewicz}.

Let $U\subset\cc^{n}$ be an open and connected neighbourhood of $0\in\cc^{n}$;
let $f=(f_{0},\ldots,f_{m})\ne0$ be a sequence of holomorphic functions
on $U$. Put $S=\{x\in U:f(x)=0\}$ and
\[
E(f)=\{(x,u)\in U\times\pm:f_{i}(x)u_{j}=f_{j}(x)u_{i},i,j=0,\ldots,m\},
\]
 where $u=[u_{0}:\cdots:u_{m}]\in\pm$.

Let $Y$ be the closure of $E(f)\setminus S$ in $U\times\pm$. The
natural projection 
\[
\pi:Y\to U
\]
is called the \emph{(elementary) blowing-up of $U$ by means of 
 $f_{0},\ldots,f_{m}$.} The ana\-lytic subset $S$ is called
a \emph{centre of the blowing-up }and its inverse image $\pi^{-1}(S)\subset Y$
is called the \emph{exceptional set }of the blowing-up.
\begin{prop}\label{blow}
Under above notations we have:
\begin{enumerate}
\item $Y$ is an analytic subset of $U\times\pm$;\label{enu:-is-an}
\item $\pi$ is proper, its range is $U$ and the restriction $\pi_{|Y\setminus\pi^{-1}(S)}$
is a biholomorphism onto $U\setminus S$;\label{enu:-is-proper,}
\item $Y$ is irreducible;\label{enu:-is-irreducible}
\item The exceptional set $\pi^{-1}(S)$ is analytic in $U\times\pm$ and
it is of pure dimension $n-1$.\label{enu:The-exceptional-set}
\end{enumerate}
\end{prop}
\begin{proof}
Although the above proposition is well known, we think that point
\eqref{enu:The-exceptional-set} is worth proving. Let us consider
the analytic map
\[
F\colon U\times\pm\ni(x,u)\mapsto(f(x),u)\in\cc^{m+1}\times\pm.
\]
 Let $y_{0},\ldots,y_{m}$ be coordinates in $\cc^{m+1}$. If we denote
by $\pi_{m+1}\colon\Pi_{m+1}\to\cc^{m+1}$ the blowing-up of $\cc^{m+1}$
by means of $y_{0},\ldots,y_{m}$ then for the restriction $\widetilde{f}=F_{|Y}$
we get the following commutative diagram of analytic maps:

\begin{center}$\begin{CD}Y @>\widetilde{f}>> \Pi_{m+1} \\
 @VV\pi V @VV\pi_{m+1} V\\
U @>f>> \cc^{m+1}
\end{CD}$\end{center}

Take $(x_{0},u_{0})\in\pi^{-1}(0)$. Let $\Omega\subset\Pi_{m+1}$
be a neighbourhood of $(0,u_{0})$, $h\colon\Omega\to\cc$ an analytic
function such that 
\[
\pi_{m+1}^{-1}(0)\cap\Omega=\{(y,u)\in\Omega:h(y,u)=0\}.
\]

Let $\widetilde{\Omega}\subset Y$ be a neighbourhood of $(x_{0},u_{0})$
such that $\widetilde{f}(\widetilde{\Omega})\subset\Omega$. Since
$\widetilde{f}^{-1}(\pi_{m+1}^{-1}(0))=\pi^{-1}(S)$ we get
\[
\pi^{-1}(S)\cap\widetilde{\Omega}=\{(x,u)\in\widetilde{\Omega}:h\circ\widetilde{f}(x,u)=0\}.
\]

Thus there exists a neighbourhood $\Delta\subset U\times\pm$ of $(x_{0},u_{0})$
and an analytic set $V\subset\Delta$ of pure dimension $n+m-1$ such
that 
\[
\pi^{-1}(S)\cap\Delta=V\cap Y\cap\Delta.
\]

This gives
\[
\dim_{(x_{0},u_{0})}\pi^{-1}(S)\geqslant\dim_{(x_{0},u_{0})}Y-1=n-1.
\]

Since $Y$ is irreducible and $\pi^{-1}(S)\varsubsetneq Y$ we get
that $\dim_{p}\pi^{-1}(S)=n-1$ for any $p\in\pi^{-1}(S)$. This ends
the proof.
\end{proof}

\section{Proof of Theorem \ref{thm:Main}}
\begin{lem}
\label{lem:blow lemma}Let $F\colon(\cc^{n}\times\cc,0)\to(\cc^{m+1},0)$,
$m\geqslant n$ be a holomorphic map. Assume that $0$ is an isolated
point of $F_{t}^{-1}(0)$ for $|t|<\delta$. Then there exists $\delta>\epsilon>0$
and a complex line $V\subset\cc^{m+1}$, such that $V\cap C_{F_{t}}=\{0\}$
for $|t|<\epsilon$.\end{lem}
\begin{proof}
Let $F\colon U\to\cc^{m+1}$, where $U\subset\cc^{n}\times\cc$ is
a connected neighbourhood of the origin. Put $S=\{(z,t)\in U:F(z,t)=0\}$
and let $\pi\colon U\times\pm\supset Y\to U$ be the elementary blowing-up
of $U$ by $F$. By Proposition \ref{blow} its exceptional set $E:=\pi^{-1}(S)$ is an
analytic set of pure dimension $n$. Let $\mathcal{E}$ be a set of
those irreducible components $W$ of $E$ for which origin in $\cc^{n+1}$
is an accumulation point of $\pi(W)\cap(\{0\}\times\cc)$. Then $\mathcal{E}$
is finite. Denote by $\widetilde{C_{F_{t}}}$ the image of the cone
$C_{F_{t}}$ in $\pm$. Observe that
\[
\{(0,t)\}\times\widetilde{C_{F_{t}}}\subset\bigcup\mathcal{E},\quad|t|<\delta.
\]
On the other hand for any $W\in\mathcal{E}$ we have 
\[
\dim W\cap(\{0\}\times\pm)\leqslant n-1<m.
\]
 Thus there exists $\epsilon>0$ and an open set $G\subset\pm$ such
that 
\[
(\{(0,t)\}\times G)\cap\bigcup\mathcal{E}=\emptyset,\quad0<|t|<\epsilon
\]
As a result if $V$ is a line in $\cc^{m+1}$ corresponding to some
point in $G$ then $V\cap C_{F_{t}}=\{0\}$ for $0<|t|<\epsilon$.
Since $G$ is not a subset of $C_{F_{0}}$ we get the assertion. 
\end{proof}

\begin{proof}[Proof of Theorem \ref{thm:Main}]
Induction on $m$. In the case $m=n$ there is nothing to prove.
Let us assume that the assertion is true for some $m\geqslant n$
and let $F\colon(\cc^{n}\times\cc,0)\to(\cc^{m+1},0)$ be a holomorphic
map such that the ideals $\left(F_{t}\right)\mathcal{O}_{n}$ are
$\mathfrak{m}_{n}$-primary. By Lemma \ref{lem:blow lemma} there
exists $\epsilon>0$ and a linear mapping $\pi'\colon\cc^{m+1}\to\cc^{m}$
such that $\ker\pi'\cap C_{F_{t}}=\left\{ 0\right\} $ for $|t|<\epsilon$.
Thus, by Corollary \ref{cor:cone-reduction} the ideal $\left(\pi'\circ F_{t}\right)\mathcal{O}_{n}$
is a reduction of $\left(F_{t}\right)\mathcal{O}_{n}$. On the other
hand, by induction hypothesis, there exists a linear map $\pi''\colon\mathbb{C}^{m}\to\mathbb{C}^{n}$
such that $\left(\pi''\circ\pi'\circ F_{t}\right)\mathcal{O}_{n}$
is a reduction of $\left(\pi'\circ F_{t}\right)\mathcal{O}_{n}$ for
small $t$. Thus if we put $\pi:=\pi''\circ\pi'$ we get the assertion.
\end{proof}
\bibliographystyle{plain}
\bibliography{/home/tomek/documents/tex/rezultaty/bibliografia.bib}

\bigskip{}

\end{document}